\documentclass[11pt]{llncs}
\usepackage[utf8]{inputenc}
\usepackage[T1]{fontenc}
\usepackage{booktabs}
\usepackage{cite}
\usepackage{dingbat}
\usepackage[normalem]{ulem}
\usepackage{mathtools,amssymb,amsmath}
\usepackage{enumerate}
\usepackage{varioref}
\usepackage[makeroom]{cancel}
\usepackage[pdfpagelabels=true]{hyperref}
\usepackage{cleveref}
\usepackage{xcolor}
\usepackage{algpseudocode}
\usepackage{algorithm}
\usepackage{listings}
\usepackage{csquotes}
\usepackage{url}
\usepackage{caption}
\usepackage{subcaption}
\usepackage{pifont}% http://ctan.org/pkg/pifont
\newcommand{\cmark}{\ding{51}}%
\newcommand{\xmark}{\ding{55}}%

\hypersetup{
	linktoc = page,
	pdfpagemode = UseNone,
	colorlinks,
	linkcolor={red!50!black},
	citecolor={blue!50!black},
	urlcolor={blue!80!black}
}

\usepackage{paralist}
\usepackage{ragged2e}
\usepackage{tikz}
\usepackage{pgfplots}
\usepackage{tikz-qtree}
\usepackage[linguistics]{forest}
\usetikzlibrary{arrows}
\usetikzlibrary{backgrounds,decorations.pathreplacing, bending}
\usepackage{breqn}
\definecolor{identifiercolor}{rgb}{.4,.6,.56}
\definecolor{stringcolor}{gray}{0.5}
\definecolor{inactivecolor}{rgb}{0.15,0.15,0.5}
\usepackage{mathtools}
\newcommand\cplus{\mathbin{\raisebox{-\height}{$+$}}}
\newcommand\contdots{\raisebox{-\height}{$\vphantom{+}\dotsm$}}

\usepackage{listings}
\begin{document}

    \title{Pattern Recognition Experiments on Mathematical Expressions}

    \author{David Naccache\inst{1} \and Ofer Yifrach-Stav\inst{1}} 

    \institute{
    DI\'ENS, \'ENS, CNRS, PSL University, Paris, France\\
45 rue d'Ulm, 75230, Paris \textsc{cedex} 05, France\\
\email{\url{ofer.friedman@ens.fr}}, 
\email{\url{david.naccache@ens.fr}}
}

\maketitle         

\begin{abstract}

We provide the results of pattern recognition experiments on mathematical expressions. 

We give a few examples of conjectured results. None of which was thoroughly checked for novelty. We did not attempt to prove all the relations found and focused on their generation.

\end{abstract}
 
\section{Introduction}

Pattern recognition is a process that involves identifying rules in data and matching them with particular case information. Pattern recognition can be seen as a type of machine learning, as it uses machine learning algorithms to recognize patterns in data. This process is characterized by the ability to learn from data, recognize familiar patterns, and recognize patterns even if they are partially visible.

Very schematically, there are three main types of pattern recognition heuristics: statistical pattern recognition, syntactic pattern recognition, and neural pattern recognition.
\begin{itemize}
    \item \textsl{Statistical pattern recognition} involves using particular case data to learn from examples and generalize rules to new observations.
    \item \textsl{Syntactic pattern recognition} (a.k.a structural pattern recognition), involves identifying patterns based on simpler sub-patterns called primitives. For example, opcodes can be seen as primitives that connect to form programs.
    \item \textsl{Neural pattern recognition} relies on artificial neural networks, which are made up of many simple processors and their connections. These networks can learn complex nonlinear input-output relationships and adapt to data through sequential training procedures.
\end{itemize}

Most pattern recognition heuristics proceed by two steps:
\begin{itemize}
    \item An \textsl{Explorative Stage} that seeks to identify patterns
    \item A \textsl{Descriptive Stage} that categorizes patterns found during exploration
\end{itemize}

In this work we provide the results of the explorative stage of syntactic pattern recognition on mathematical expressions. Given the nature of the objects we work on (conjectures) the descriptive stage is left to a human.

We give a few examples of conjectured results. None of which was thoroughly checked for novelty. We did not attempt to prove all the relations found and focused on their generation.

\section{The Pattern Recognition Algorithm}

The pattern recognition algorithm has two components called the \textsl{generalizer} and the \textsl{identifier}. 

The generalizer departs from a known continued fraction or a mathematical expression (a particular case) and automatically parameterizes parts of it. The parameterized parts are \textsl{target ingredients} tagged by the user. For each set of particular parameter values (taken over search space), approximated values of the formula are collected for later analysis. 

Target ingredients are replaced by progressions, denoted by $\mu_{\mathbf{u}}(i)$, which can be constant, (alternating) arithmetic, geometric, harmonic or exponential depending on the parameter choices. Those are captured by the general formula:

$$\mu_{\mathbf{u}}(i)=u_4i^{u_5}+(u_0+i u_1)^{u_3}u_2^i$$

For instance, the Ramanujan Machine Project \cite{rama,ref1,ref2} re-discovered an already known relation involving $e^\pi$. Namely, that the continued fraction defined by $b_n=n^2+4$ and $a_n=2 n+1$ converges to:

\begin{equation*}\label{eq:gcf}
\frac{2 \left(e^{\pi }+1\right)}{e^{\pi }-1}=1+{\cfrac {1^2+4}{3+{\cfrac {2^2+4}{5+{\cfrac {3^2+4}{7+{\cfrac {4^2+4}{9+\ddots \,}}}}}}}}
\end{equation*}

A natural tagging query of this identity for search by the user might hence be:

\begin{equation*}\label{eq:gcf}
Q(\mathbf{u})=\mu_{\mathbf{u}}(0)+{\cfrac{\mu_{\mathbf{v}}(0)}{\mu_{\mathbf{u}}(1)+{\cfrac {\mu_{\mathbf{v}}(1)}{\mu_{\mathbf{u}}(2)+{\cfrac {\mu_{\mathbf{v}}(2)}{\mu_{\mathbf{u}}(3)+{\cfrac {\mu_{\mathbf{v}}(3)}{\mu_{\mathbf{u}}(4)+\ddots \,}}}}}}}}
\end{equation*}

With $$\mathbf{u}=\{\mathbb{Q},\mathbb{Q},1,1,0,0\}\mbox{~~and~~} \mathbf{v}=\{\mathbb{Z},0,1,1,\mathbb{Q},\mathbb{N}\}$$

That is:

$$\mu_{\mathbf{u}}(i)=(\mathbb{Q}+i\mathbb{Q}) \mbox{~~and~~} \mu_{\mathbf{v}}(i)=\mathbb{Q}i^{\mathbb{N}}+\mathbb{Z}$$

When this is done, the program varies the progressions' parameters over the chosen search spaces and collects sequences of resulting values. The tests that we list here are of course non limitative and many other variants can be added to the proposed heuristic.

\remark{Obviously, we are quickly limited by the increasing complexity due to nested loops running over the parameters of the expressions (i.e. the $u_i$s).}\smallskip

\remark{At the risk of overlooking some gold nuggets, when we explore $\mathbb{Q}$ we start by exploring $\mathbb{N}$ and if the search is conclusive, we refine it by increments of $1/6$ which have the advantage of exploring units, halves and thirds at the cost of a small multiplicative factor of $6$. If interesting results are found with increments of 6 the step is refined to $1/30$ and to Farey sequences.}\smallskip

The sequences obtained by varying those parameters are fed into the identifier for possible recognition. To detect conjectures the identifier performs a number of tests on the obtained sequences. Tests belong to two categories: \textsl{morphological tests} and \textsl{serial tests}. Morphological tests are applied to very few individual results and try to spot their characteristics. Serial tests are applied to more results and seek to discover relationships between them.\smallskip

\textbf{Algebraic number identification (ANI)}: Collect 10 convergence limits $Q_0,Q_1,\ldots Q_9$ and, using LLL \cite{LLL}, check if any of those $Q_i$s is the root of a small degree ($\leq 10$) polynomial. If so, check that RNI failed before returning \textsf{true} to avoid multiple alerts as rationals are also algebraic. This is a morphological test. The degree $10$ was chosen arbitrarily and can be changed at wish (provided that the precision is matched to the degree).\smallskip

\textbf{Rational number identification (RNI)}: Collect 10 convergence limits $Q_0,Q_1,\ldots Q_9$ and, using LLL, check if any of those $Q_i$s is a good approximation of a rational number having a (abnormally) small numerator and a small denominator. This is a morphological test.\smallskip

\textbf{Constant presence identification (CPI)}: Collect 10 convergence limits $Q_0,Q_1,\ldots Q_9$. Consider the 45 pairs $P_1,P_2$ formed from those $Q_i$s. Using LLL, check the assumption that there is at least one pair of the form:

$$
P_1=\frac{a_1+ b_1 U}{c_1+ d_1 U}\mbox{~~and~~}P_2=\frac{a_2+ b_2 U }{c_2+ d_2 U}
$$

Where $U\not\in\mathbb{Q}$ and $a_1,b_1,c_1,d_1,a_2,b_2,c_2,d_2\in \mathbb{Z}$.

Solving for $U$ and equating we get:

$$
a_2 b_1 - a_1 b_2 + (b_2 c_1 - a_2 d_1) P_1  + (a_1 d_2 - b_1 c_2)P_2+ 
 (c_2 d_1  - c_1 d_2) P_1 P_2=0$$
 
Hence, when called with on input $1,P_1,P_2,P_1 P_2$ LLL will return an abnormally short vector if the coefficients are small (as is usually the case in remarkable identities). This is a morphological test.\smallskip

\textbf{Constant to exponent identification (CEI)}: Collect 10 convergence limits $Q_0,Q_1,\ldots Q_9$. Consider the 7 quadruples $P_1,P_2,P_3,P_4$ formed by successive $Q_i$s\footnote{namely: $\{0,1,2,3\}$,$\{1,2,3,4\}$,$\{2,3,4,5\}$,$\{3,4,5,6\}$,$\{4,5,6,7\}$,$\{5,6,7,8\}$,$\{6,7,8,9\}$}. 

Here we assume that at successive ranks the limits are of the form:

$$
P_k=\frac{a_k+ b_k U^k}{c_k+ d_k U^k}$$

Which implies that:

$$U^k = \frac{a_k-c_k P_k}{d_k P_k-b_k}
$$

If follows that:

$$U= \frac{(a_{k+1}-c_{k+1} P_{k+1})({d_k P_k-b_k})}{(d_{k+1} Q_{k+1}-b_{k+1})({a_k-c_k Q_k})}$$

$$\frac{(a_{k+3}-c_{k+3} P_{k+3})({d_{k+2} P_{k+2}-b_{k+2}})}{(d_{k+3} P_{k+3}-b_{k+3})({a_{k+2}-c_{k+2} P_{k+2}})} = \frac{(a_{k+1}-c_{k+1} P_{k+1})({d_k P_k-b_k})}{(d_{k+1} P_{k+1}-b_{k+1})({a_k-c_k P_k})}$$

Let:

$$S_1=\{P_k,P_{k+1},P_{k+2},P_{k+3}\}$$
$$S_2=\{P_k P_{k+1},P_k P_{k+2},P_{k+1} P_{k+2},P_k P_{k+3},P_{k+1} P_{k+3},P_{k+2} P_{k+3}\}$$
$$S_3=\{P_k P_{k+1} P_{k+2},P_k P_{k+1} P_{k+3},P_k P_{k+2} P_{k+3},P_{k+1} P_{k+2} P_{k+3}\}$$
$$S=S_1 \cup S_2 \cup S_3 \cup \{1,P_k P_{k+1} P_{k+2} P_{k+3}\}$$

When called with on input $S$ LLL will return an abnormally short vector (as is usually the case in remarkable identities). This is a morphological test.\smallskip

\remark{Both CPI and CEI can be generalized to detect the presence of multiple unknown constants in an expression (i.e. $U_1,U_2,\ldots$) or even the presence of common constants in different continued fractions. We did not implement this generalization. Following those tests we can compute a numerical approximation of $U$ and attempt to look it up\footnote{e.g. on \url{https://wayback.cecm.sfu.ca/projects/ISC/ISCmain.html}}}.\smallskip

\textbf{Known constant identification (KCI)}: Let $L$ be the following set of usual constants:

$$L=\{1,\sqrt{\pi},\pi,\pi^2,\pi^3,\zeta(3),\zeta(5),\zeta(7),\sqrt{e},e,e^2,e^3,\phi^2,\gamma,G,\ln{2},\ln{3},\ln{5}\}$$

Collect 10 convergence limits $Q_0,Q_1,\ldots Q_9$. Check using LLL if any of the $Q_i$ is a number of the form:

$$
Q_i\sum_{j} a_j L_j=\sum_{j}  b_j L_i \mbox{~~for~~}a_1,a_2,\ldots,b_1,b_2\ldots\in\mathbb{Z}
$$

If the solution only involves $1$, a \textsf{false} is returned. Note that as $L$ increases the required precision must also be increased to prevent spotting artefacts. In practice we (manually) select only a subset of $L$ before running the KCI test according to the nature of the constants appearing the in the particular case. Note that KCI and CPI can have overlapping responses.\smallskip

\textbf{Rational fraction progression (RFP)}: In this test we seek to see if when all $u_i$ except one (say $\bar{u}$) are kept constant, the continued fraction's limit $Q(\bar{u})$ is a ratio of two polynomials in $\bar{u}$ with integer coefficients. This is done by a non linear model fit. The fit residuals serve as a measure of the verdict's likelihood. This is a serial test.\smallskip

\textbf{Exponential function progression (EFP)}: In this test we seek to see if when all $u_i$ except one (say $\bar{u}$) are kept constant, the continued fraction's limit $Q(\bar{u})$ is a function of the form $b a^{\bar{u}}$ with rational coefficients. This is done by a non linear model fit and rationality detection on $a,b$. The fit residuals serve as a measure of the verdict's likelihood. If $ab=0$ return \textsf{false} to avoid reporting the same result as the RFP. This is a serial test.\smallskip

\textbf{Inverse exponential progression (IEP)}: In this test we seek to see if when all $u_i$ except one (say $\bar{u}$) are kept constant, the continued fraction's limit $Q(\bar{u})$ is a function of the form $b a^{1/\bar{u}}$ with rational coefficients. This is done by a non linear model fit and rationality detection on $a,b$. The fit residuals serve as a measure of the verdict's likelihood. If $ab=0$ return \textsf{false} to avoid reporting the same result as the RFP. This is a serial test.\smallskip

\textbf{Power plus constant progression (PCP)}: In this test we seek to see if when all $u_i$ except one (say $\bar{u}$) are kept constant, the continued fraction's limit $Q(\bar{u})$ is a function of the form $b \bar{u}^a+c$ with rational coefficients. This is done by a non linear model fit and rationality detection on $a,b,c$. The fit residuals serve as a measure of the verdict's likelihood. If $b=0$ return \textsf{false} to avoid reporting the same result as the RFP. This is a serial test.\smallskip

\textbf{Root plus constant progression (RCP)}:  In this test we seek to see if when all $u_i$ except one (say $\bar{u}$) are kept constant, the continued fraction's limit $Q(\bar{u})$ is a function of the form $b \sqrt[a]{\bar{u}}+c$ with rational coefficients. This is done by a non linear model fit and rationality detection on $a,b,c$. The fit residuals serve as a measure of the verdict's likelihood. If $ab=0$ return \textsf{false} to avoid reporting the same result as the RFP. This is a serial test.\smallskip

\section{Continued Fractions Converging to $2 u (e^{u\pi }+1)/(e^{u\pi }-1)$}

It appears that the relation:

\begin{equation*}\label{eq:gcf}
\frac{2 \left(e^{\pi }+1\right)}{e^{\pi }-1}=1+{\cfrac {1^2+4}{3+{\cfrac {2^2+4}{5+{\cfrac {3^2+4}{7+{\cfrac {4^2+4}{9+\ddots \,}}}}}}}}
\end{equation*}

is the first in an infinite family:

\begin{equation*}\label{eq:gcf}
\frac{2 u \left(e^{u \pi }+1\right)}{e^{u \pi }-1}=1+{\cfrac {1^2+4 u^2}{3+{\cfrac {2^2+4 u^2}{5+{\cfrac {3^2+4 u^2}{7+{\cfrac {4^2+4 u^2}{9+\ddots \,}}}}}}}}
\end{equation*}

\begin{table}[]
    \centering
    \begin{tabular}{|c|c|c|c|c|c|c|c|c|c|}\hline
~~ANI~~&~~RNI~~&~~CPI~~&~~CEI~~&~~KCI~~&~~RFP~~&~~EFP~~&~~IEP~~&~~PCP~~&~~RCP~~\\\hline\hline
\xmark &\cmark &\cmark &\cmark & \xmark&\xmark &\xmark &\xmark &\xmark & \cmark\\\hline
    \end{tabular}
    \caption{Test Results}
    \label{tests1}
\end{table}

Indeed, (RCP) linear variations in $u$ cause identifiable $O(\sqrt{u})$ variations in the limit. This is because very quickly:

$$\lim_{u\rightarrow \infty}\frac{e^{u\pi }+1}{e^{u\pi }-1}= 1$$

This has the somewhat adverse effect of making the RNI positive very quickly as well.

The final form is detected thanks to the CEI test.

\subsubsection{By-product: } Because this holds for $u\in\mathbb{C}^*$, we get a few seemingly ``mysterious'' corollary identities such as:

\begin{equation*}\label{eq:gcf}
\frac{2\left(e+1\right)}{\pi (e-1)}=1+{\cfrac {1^2+4/\pi^2}{3+{\cfrac {2^2+4/\pi^2}{5+{\cfrac {3^2+4/\pi^2}{7+{\cfrac {4^2+4/\pi^2}{9+\ddots \,}}}}}}}}
\end{equation*}

\begin{equation*}\label{eq:gcf}
\frac{6\ln 2}{\pi}=1+{\cfrac {1^2+4\ln^2 2/\pi^2}{3+{\cfrac {2^2+4\ln^2 2/\pi^2}{5+{\cfrac {3^2+4\ln^2 2/\pi^2}{7+{\cfrac {4^2+4\ln^2 2/\pi^2}{9+\ddots \,}}}}}}}}
\end{equation*}

\subsubsection{Implementation}
~\\

\begin{lstlisting}[extendedchars=true,language=Mathematica]
f[x_, {m_, d_}] := m/(d + x);
For[t = 0, t <= 5,
  den = Table[2 n + 1, {n, 1, 20000}];
  num = Table[n^2 + (2 t)^2, {n, 1, 20000}];
  r = 1 + (Fold[f, Last@num/Last@den, Reverse@Most@Transpose@{num, den}]);
  e = 2 t (1 + (E^Pi)^t)/((E^Pi)^t - 1);
  Print[{e, 2 n + 1, n^2 + (2 t)^2, N[{r, e}, 20]}];
  t += 1/2];
\end{lstlisting}

\section{Continued Fractions Converging to Polynomial Roots}

It is very well known that: 

\[\frac{\sqrt{5}-1}{2}=\frac{1}{1}\cplus\frac{1}{1}\cplus\frac{1}{1}\cplus\frac{1}{1}\cplus\frac{1}{1}\cplus
\frac{1}{1}\cplus\frac{1}{1}\cplus\frac{1}{1}\cplus
\contdots\]

We tag\footnote{Adding a $1+$ by commodity which does not change anything about the infinite convergence.}:

\[Q(\mathbf{u})=1+\frac{\mu_{\mathbf{u}}(0)}{\mu_{\mathbf{u}}(0)}\cplus\frac{\mu_{\mathbf{u}}(0)}{\mu_{\mathbf{u}}(0)}\cplus\frac{\mu_{\mathbf{u}}(0)}{\mu_{\mathbf{u}}(0)}\cplus\frac{\mu_{\mathbf{u}}(0)}{\mu_{\mathbf{u}}(0)}\cplus\frac{\mu_{\mathbf{u}}(0)}{\mu_{\mathbf{u}}(0)}\cplus
\frac{\mu_{\mathbf{u}}(0)}{\mu_{\mathbf{u}}(0)}\cplus\frac{\mu_{\mathbf{u}}(0)}{\mu_{\mathbf{u}}(0)}\cplus\contdots\]

With:

$$\mathbf{u}=\{\mathbb{Q},0,1,1,0,0\} \Rightarrow \mu_{\mathbf{u}}(i)=\mathbb{Q}$$

It appears that for $u\in\mathbb{Q}/[-4,0]$ LLL identifies that the limit is a root of a second degree polynomial, namely:

\[Q(u)=1+\frac{u}{u}\cplus\frac{u}{u}\cplus\frac{u}{u}\cplus\frac{u}{u}\cplus\frac{u}{u}\cplus\frac{u}{u}\cplus\frac{u}{u}\cplus\contdots\]

$$Q(u)^2+u(Q(u)-1)=0$$

Which is trivial to prove by pushing the $u$ into the continued fraction.

\begin{table}[]
    \centering
    \begin{tabular}{|c|c|c|c|c|c|c|c|c|c|}\hline
~~ANI~~&~~RNI~~&~~CPI~~&~~CEI~~&~~KCI~~&~~RFP~~&~~EFP~~&~~IEP~~&~~PCP~~&~~RCP~~\\\hline\hline
\cmark &\xmark &\cmark &\cmark & \xmark&\xmark &\xmark & \xmark&\xmark & \xmark\\\hline
    \end{tabular}
    \caption{Test Results}
    \label{tests2}
\end{table}

The CPI is positive because for $u=1$ and $u=5$ the respective values of $Q(u)$ comprise the common value $\sqrt{5}$.

\subsubsection{Implementation}
~\\

\begin{lstlisting}[extendedchars=true,language=Mathematica]
f[x_, {m_, d_}] := m/(d + x);
For[L = -20, L <= 20 ,
  If[-4 <= L <= 0, L = 2/3]; 
  num = den = Table[L, {n, 1, 200}];
  r = Fold[f, Last@num/Last@den, Reverse@Most@Transpose@{num, den}];
  Print[{L, N[r^2 + L (r - 1)]}];
  L += 2/3];
\end{lstlisting}

\section{Continued Fractions Converging to $e^{2/\kappa}$}

The following relations are well-known\footnote{\url{https://link.springer.com/content/pdf/bbm:978-94-91216-37-4/1.pdf}}:

\[e=2+\frac{1}{1}\cplus\frac{1}{2}\cplus\frac{1}{1}\cplus\frac{1}{1}\cplus\frac{1}{4}\cplus
\frac{1}{1}\cplus\frac{1}{1}\cplus\frac{1}{6}\cplus
\contdots\]

\[\sqrt{e}=1+\frac{1}{1}\cplus\frac{1}{1}\cplus\frac{1}{5}\cplus\frac{1}{1}\cplus\frac{1}{1}\cplus\frac{1}{9}\cplus
\frac{1}{1}\cplus\frac{1}{1}\cplus\frac{1}{13}\cplus
\contdots\]

\[\sqrt[3]{e}=1+\frac{1}{2}\cplus\frac{1}{1}\cplus\frac{1}{1}\cplus\frac{1}{8}\cplus
\frac{1}{1}\cplus\frac{1}{1}\cplus\frac{1}{14}\cplus\frac{1}{1}\cplus\frac{1}{1}\cplus\frac{1}{20}\cplus
\contdots\]

We hence tag the ones as constants, the progression as arithmetic and let the algorithm monitor the evolution of the limits.\smallskip

Let $b_n=1$. Define $\mu(u)=\kappa(u+1/2)-1$ for $\kappa\in\mathbb{R}$ and:

\[
a_n =
 \begin{cases}
  \mu(n/3)= \frac{\kappa(2 n+3)}{6} -1& \mbox{~if~} n \bmod 3\equiv0 \\
  1          &  \mbox{~otherwise.~}
 \end{cases}
\]

In other words, $a_n$ is the sequence:

$$a_n=\{\mu(0),1,1,\mu(1),1,1,\mu(2),1,1,\mu(3),1,1,\mu(4),1,1,\cdots\}$$

Then we detect that the continued fraction generated by $a_n,b_n$ converges to $e^{2/\kappa}$.

\begin{table}[]
    \centering
    \begin{tabular}{|c|c|c|c|c|c|c|c|c|c|}\hline
~~ANI~~&~~RNI~~&~~CPI~~&~~CEI~~&~~KCI~~&~~RFP~~&~~EFP~~&~~IEP~~&~~PCP~~&~~RCP~~\\\hline\hline
\xmark &\xmark &\xmark &\cmark & \xmark&\xmark &\xmark & \cmark&\xmark & \xmark\\\hline
    \end{tabular}
    \caption{Test Results}
    \label{tests3}
\end{table}

The CEI is positive because, for instance $(e^{2/\kappa})^2=e^{2/\kappa'}$ implies that $2/\kappa=\kappa'$ which is satisfied for several pairs of integer values.

\subsubsection{Implementation}
~\\

\begin{lstlisting}[extendedchars=true,language=Mathematica]
f[x_, {m_, d_}] := m/(d + x);
For[k = -10, k <= 10,
  phi = Table[k n + k/2 - 1, {n, 0, 2000 - 1}] ;
  num = Table[1, {n, 1, 2000}];
  den = Take[
    Flatten[Table[{phi[[i]], {1, 1}}, {i, 1, Floor[2000/3] + 1}]], {1,
      2000}];
  r = 1 + (Fold[f, Last@num/Last@den, Reverse@Most@Transpose@{num, den}]);
  v = E^(2/k);
  Print[{k, v, N[{r, v}, 20]}];
  k += 1/2];
\end{lstlisting}

\section{Continued Fractions Involving Catalan's Constant}

It is well known that:

\[2G=2-\frac{1^2}{3}\cplus\frac{2^2}{1}\cplus\frac{2^2}{3}\cplus\frac{4^2}{1}\cplus\frac{4^2}{3}\cplus
\frac{6^2}{1}\cplus\frac{6^2}{3}\cplus\frac{8^2}{1}\cplus\frac{8^2}{3}\cplus
\contdots\]

We define:

\[\Delta(u,v)=\frac{1}{2v}\times\left(\frac{1^2}{u}\cplus\frac{2^2}{v}\cplus\frac{2^2}{u}\cplus\frac{4^2}{v}\cplus\frac{4^2}{u}\cplus
\frac{6^2}{v}\cplus\frac{6^2}{u}\cplus\frac{8^2}{v}\cplus\frac{8^2}{u}\cplus
\contdots\right)\]

For all the following we observe that $\Delta(u,v)=\Delta(v,u)$.

\subsection{For $u=1$}

An exploration for $\mathbf{u}=\{0,\mathbb{N},\mathbb{N},\mathbb{N},\mathbb{Z},0\}$ reveals that for $u_0=0,u_1=2,u_2=1,u_3=2,u_4=-1,u_5=0$ we get identities when $v=4i^2-1$ with the convergence values given in Table \ref{tab:my_label}:

\begin{table}[]
    \centering
    \begin{tabular}{|c|c|c|c|}\hline
~~~$u$~~~&~~~$i$~~~&~~~$v=4i^2-1$~~~&~~~$\Delta(u,4i^2-1)=\Delta(1,4i^2-1)$~~~\\\hline\hline
1   &   0   &    -1       &$1-G$\\\hline
1   &   1   &   3       &$-8/9+G$\\\hline
1   &   2   &   15       &$209/225-G$\\\hline
1   &   3   &   35       &$-10016/11025+G$\\\hline
1   &   4   &   63       &$91369/99225-G$\\\hline
1   &   5   &  99       &$-10956424/12006225+G$\\\hline
1   &   6   &  143       &~~~$1863641881/2029052025-G$~~~\\\hline
    \end{tabular}
    \caption{The first convergence values for $u=1$}
    \label{tab:my_label}
\end{table}

Where the general formula for $i>1$ is:

$$
\Delta(1,4i^2-1) = (-1)^{i+1} \left(\sum_{k=0}^{i-1} \frac{(-1)^k}{(2k + 1)^2}-G\right)
$$

\begin{table}[]
    \centering
    \begin{tabular}{|c|c|c|c|c|c|c|c|c|c|}\hline
~~ANI~~&~~RNI~~&~~CPI~~&~~CEI~~&~~KCI~~&~~RFP~~&~~EFP~~&~~IEP~~&~~PCP~~&~~RCP~~\\\hline\hline
\xmark &\xmark &\cmark &\xmark & \cmark&\xmark &\xmark & \xmark&\xmark & \xmark\\\hline
    \end{tabular}
    \caption{Test Results}
    \label{tests4}
\end{table}

\subsubsection{Implementation}
~\\

\begin{lstlisting}[extendedchars=true,language=Mathematica]
f[x_, {m_, d_}] := m/(d + x);
For[i = 1, i < 40,
  {u, v} = {1, 4 i^2 - 1};
  num = Take[
    Prepend[Flatten[Table[{(2 n)^2, (2 n )^2}, {n, 1, 100000}]], 1] , 
    100000];
  den = Flatten[Table[{u, v}, {n, 1, 100000/2}]];
  r = Fold[f, Last@num/Last@den, Reverse@Most@Transpose@{num, den}]/2/v;
  val = (-1)^(i + 1) (Sum[(-1)^k/(2 k + 1)^2, {k, 0, i - 1}] - Catalan);
  Print[{i, v, val, N[{val, r}, 30]}];
  i++];
\end{lstlisting}

\remark{Note that the denominators of the numbers:
$$\eta(i)=\sum_{k=0}^{i-1} \frac{(-1)^k}{(2k + 1)^2}$$
are interesting by their own right. At first sight they might seem perfect squares but in reality some may contain very small prime factors to an odd power.}

\subsection{For $u=3$}

The exploration in this section is interesting. It was done manually but we would have never had the idea to probe in that specific direction without the insight for the case $u=1$ produced in the previous section.

\begin{table}[]
    \centering
    \begin{tabular}{|c|c|c|c|}\hline
~~~$u$~~~&~~~$i$~~~&~~~$f(i)$~~~&~~~$\Delta(3,f(i))$~~~\\\hline\hline
3       &   0   &    1       &$\Delta(1,\kern 0.15em-1)$\\\hline
3       &   1   &   5       &$\Delta(1,\phantom{00}3)$\\\hline
3       &   2   &   21       &$\Delta(1,\phantom{0}35)$\\\hline
3       &   3   &   33       &$\Delta(1,\phantom{0}63)$\\\hline
3       &   4   &   65       & $\Delta(1,143)$\\\hline
3       &   5   &   85       &$\Delta(1,255)$\\\hline
    \end{tabular}
    \caption{The first convergence values for $u=3$}
    \label{tab:case3}
\end{table}

The sequence $f(i)$ is \textcolor{red}{nearly} the absolute value of the OEIS sequence A006309\footnote{\url{https://oeis.org/A006309.}}:

\begin{equation*}
\begin{split}
&1, 5, 21, 33, 65, 85, 133, 161, 261, 341, 481, 533, 645,705, 901, \textcolor{red}{\xcancel{12803}}, 1281,  \\
&1541, 1633, 1825, \textcolor{red}{\xcancel{14615}}, \textcolor{red}{\xcancel{11537}}, 2581, 3201, 3333\ldots 
\end{split}
\end{equation*}

An unexplained phenomenon occurs for the ``abnormally larger'' OEIS sequence A006309 values 12803, 14615, 11537 that remains unmatched by any $\eta(i)$ value. We do not have an explanation for this phenomenon that requires further research.

\subsubsection{Implementation}
~\\

The following implementation was purposely left unoptimized for the sake of clarity. We start by generating the target values for $u=3$ and store them in an array. Then we re-generate the values for $u=1$ and match the array's contents.

\begin{lstlisting}[extendedchars=true,language=Mathematica]
AbsA006309 = 
  Abs[{1, 5, -21, 33, -65, 85, -133, 161, 261, -341, -481, 533, -645, 
    705, 901, -12803, -1281, -1541, 1633, -1825}];
t = {};
f[x_, {m_, d_}] := m/(d + x);
For[i = 1, i <= Length[AbsA006309],
  {u, v} = {3, AbsA006309[[i]]};
  num = Take[
    Prepend[Flatten[Table[{(2 n)^2, (2 n)^2}, {n, 1, 40000}]], 1], 
    40000];
  den = Flatten[Table[{u, v}, {n, 1, 40000/2}]];
  r = Fold[f, Last@num/Last@den, Reverse@Most@Transpose@{num, den}]/2/
    v;
  AppendTo[t, {AbsA006309[[i]], N[r, 30]}];
  i++];

For[j = 1, j <= Length[AbsA006309],
  If[t[[j, 1]] == 12803, 
   Print["Exception, the value 12803 is skipped."],
   For[i = 1, i <= 1000000,
    {u, v} = {1, 4 i^2 - 1};
    den = Flatten[Table[{u, v}, {n, 1, 40000/2}]];
    r = Fold[f, Last@num/Last@den, Reverse@Most@Transpose@{num, den}]/
       2/v;
    val = (-1)^(i + 1) (Sum[(-1)^k/(2 k + 1)^2, {k, 0, i - 1}] - 
        Catalan);
    If[Abs[t[[j, 2]] - r] < 10^(-6),
     Print[{i, N[r, 30], N[t[[j, 2]], 30]}, "Entry ", j, ": ", val, 
      " matched with Delta[3,", t[[j, 1]], "]"];
     i = Infinity];
    i++]];
  j++];
\end{lstlisting}

\subsection{Subsequent $u$ values.}

Table \ref{tab:over3} provides some additional examples for various $u,v$ combinations.

\begin{table}[]
    \centering
    \begin{tabular}{|c|c|c|}\hline
~~~$u$~~~&~~~$v$~~~&~~~$\Delta(u,v)$~~~\\\hline\hline
 5       &   7     &$\Delta(1,\phantom{1}15)$\\\hline
 5       &  39     &$\Delta(1,143)$\\\hline
 5       &  51     &$\Delta(1,255)$\\\hline 
 7       &   9     &$\Delta(1,\phantom{1}35)$\\\hline
 9       &  11     &$\Delta(1,\phantom{1}63)$\\\hline
11       &  13     &$\Delta(1,\phantom{1}99)$\\\hline
13       &  15     &$\Delta(1,143)$\\\hline
    \end{tabular}
    \caption{Other convergence values.}
    \label{tab:over3}
\end{table}

\subsection{Variations in the numerator.}

Let, for instance, $(u,v)=(1,3)$. Removing the $1/(2v)$ factor in $\Delta$ and replacing the $(2n)^2$ by $(n-i)^2$ we get convergence to:

$$1,\frac{4}{5},\frac{31}{51},\frac{16}{33},\frac{355}{883},\frac{11524}{33599},\frac{171887}{575075},\frac{10147688}{38326363},\ldots$$

With the limits being quickly reached after a constant number of terms in the continued fraction.

\subsubsection{Implementation}
~\\

\begin{lstlisting}[extendedchars=true,language=Mathematica]
For[i = 1, i < 20,
 f[x_, {m_, d_}] := m/(d + x);
 
 num = Take[
   Prepend[Flatten[Table[{(n - i)^2, (n - i)^2}, {n, 1, 400}]], 1] , 
   400];
 den = Flatten[Table[{1, 3}, {n, 1, 400/2}]];
 r = (Fold[f, Last@num/Last@den, Reverse@Most@Transpose@{num, den}]);
 Print[r];
 i += 1]
\end{lstlisting}

\section{Generalized Cloître Series}

In an unpublished note \cite{cloitre}, Benoît Cloître gives a beautiful BBP formula for $\pi^2$ based on the identity:  

$$
\sum_{k=1}^{\infty}\frac{\cos(i k \pi) \left(2 \cos(j\pi)\right)^k}{k^2}=\left(\ell \pi\right)^2
$$

Here are some $i,j,\ell$ combinations detected automatically:

\begin{table}[]
    \centering
\begin{tabular}{|c|c|c|}\hline
~~~$i/\ell$~~~&~~~$j/\ell$~~~&~~~$1/\ell$~~~\\\hline\hline
$ 11$ & $5$ & $8$ \\\hline
$ 14$ & $6$ & $10$ \\\hline
$ 23$ & $9$ & $16$ \\\hline
$ 26$ & $10$ & $18$ \\\hline
$ 22/5$ & $8/5$ & $3$ \\\hline
$ 31$ & $9$ & $20$ \\\hline
$ 28$ & $8$ & $18$ \\\hline
$ 19$ & $5$ & $12$ \\\hline
$ 16$ & $4$ & $10$ \\\hline
$ 13$ & $3$ & $8$ \\\hline
\end{tabular}
    \caption{Example relations for which $i+j=2$}
    \label{tab:sum2}
\end{table}

\begin{table}[]
    \centering
\begin{tabular}{|c|c|c|}\hline
~~~$i/\ell$~~~&~~~$j/\ell$~~~&~~~$1/\ell$~~~\\\hline\hline
$ 76$ & $16$ & $30$ \\\hline
$ 46$ & $10$ & $18$ \\\hline
$ 41$ & $9$ & $16$ \\\hline
$ 26$ & $6$ & $10$ \\\hline
$ 21$ & $5$ & $8$ \\\hline
\end{tabular}
    \caption{Example relations for which $i+j\neq 2$}
    \label{tab:sumnot2}
\end{table}

A simple rule allowing to generate many identities consists in fixing a factional step $1/u$, letting $i=\kappa/u$ for $\pi/3\leq i \leq 2\pi/3$ and calculating the limit for $\{i,j\}=\{\kappa u,2-\kappa u\}$ (e.g. Table \ref{tab:sum2}). However, limits for which $i+j\neq 2$ exist as well (e.g. Table \ref{tab:sumnot2}).

\begin{table}[]
    \centering
    \begin{tabular}{|c|c|c|c|c|c|c|c|c|c|}\hline
~~ANI~~&~~RNI~~&~~CPI~~&~~CEI~~&~~KCI~~&~~RFP~~&~~EFP~~&~~IEP~~&~~PCP~~&~~RCP~~\\\hline\hline
\xmark &\xmark &\cmark &\cmark & \cmark&\xmark &\xmark & \xmark&\xmark & \xmark\\\hline
    \end{tabular}
    \caption{Test Results}
    \label{tests5}
\end{table}

\section{Conclusion \& further research}

The results given in this paper show that pattern matching can obviously be of help in detecting new mathematical conjectures. The very basic processes described in the previous sections can be improved and generalized in a number of ways. The first is obviously an enriching of the collection of tests. The second is deeper exploration which is highly dependent on the computational capabilities at hand. Finally the interpretation of results and the early pruning of less probable branches in the potential conjecture tree can also bring efficiency and pertinence in the discovered relations. 

\bibliographystyle{abbrv}
\bibliography{biblio.bib}
\end{document}